\documentclass{amsart}

  \def\F{{\bf F}} \def\N{{\bf N}} 
 \def\Q{{\bf Q}}  \def\Z{{\bf Z}}

\def\dd#1#2{\frac{\partial #1}{\partial #2}}

\def\qed{\rule{2mm}{2mm}\bigskip}

\newtheorem*{Conjecture}{Conjecture}

\newtheorem{Lemma}{Lemma}
\newtheorem*{Proposition}{Proposition}
\newtheorem*{Theorem}{Theorem}

\begin{document}

\title{Reducibility of polynomials $f(x,y)$ modulo $p$}
\author{Wolfgang M. Ruppert}
\address{Mathematisches Institut, Universit\"at
Erlangen--N\"urnberg, Bismarckstra{\ss}e $1\frac{1}{2}$, D-91054
Erlangen, Germany}
\email{ruppert@mi.uni-erlangen.de}
\date{August 5, 1998}

\begin{abstract} We consider absolutely irreducible
polynomials $f\in\Z[x,y]$ with $\deg_xf=m$, $\deg_yf=n$ and height 
$H$. We show that 
for any prime $p$ with $p>[m(n+1)n^2+(m+1)(n-1)m^2]^{mn+\frac{n-1}{2}}
\cdot H^{2mn+n-1}$ the
reduction $f\bmod p$ is also absolutely irreducible.
Furthermore if the Bouniakowsky conjecture is true we show
that there are infinitely many absolutely irreducible polynomials
$f\in\Z[x,y]$ which are reducible mod $p$ where $p$ is a prime
with $p\ge H^{2m}$.
\end{abstract}

\maketitle

%\tableofcontents

\section{Introduction}

It is well known that for an absolutely irreducible polynomial
$f\in\Z[x,y]$ the reduction $f\bmod p$ is also absolutely irreducible
if the prime $p$ is large enough. For small $p$ the polynomial $f\bmod
p$ may be reducible. 
E.g. $f=x^9y-9x^9-2x+9y+2$ is absolutely irreducible over
$\Q$ but reducible modulo $p=186940255267545011$ where 
$x-93470127633772547$ divides $f\bmod p$. 
It is natural to ask how large $p$
has to be to be sure that $f\bmod p$ is absolutely irreducible. 
In \cite{R} we showed that 
$$p>d^{3d^2-3}\cdot H(f)^{d^2-1}$$
is sufficient for absolute irreducibility mod $p$ where $d$ is the total
degree of $f$ and $H(f)$ the height\footnote{the height of a
polynomial $f=\sum_{i,j}a_{ij}x^iy^j\in \Z[x,y]$ as we use it is
defined by $H(f)=\max_{i,j}|a_{ij}|$.} of $f$. 
Sometimes it is more natural to consider the polynomial having 
degree $m$ in $x$ and
$n$ in $y$. For this case Zannier \cite{Z} has shown that 
$$p>e^{12n^2m^2}(4n^2m)^{8n^2m}\cdot H(f)^{2(2n-1)^2m}$$
is sufficient for absolute irreducibility $\bmod p$. 
Our aim is to improve Zannier's estimate by showing the following 
theorem:

\begin{Theorem} Let $f\in\Z[x,y]$ be an absolutely irreducible
polynomial with degree $m\ge 1$ in $x$, $n\ge 1$ in $y$ and 
height $H(f)$. If $p$ is a prime with 
$$p>[m(n+1)n^2+(m+1)(n-1)m^2]^{mn+\frac{n-1}{2}}\cdot H(f)^{2mn+n-1}$$
then the reduced polynomial $f\bmod p$ is also absolutely irreducible.
\end{Theorem}

The basic ingredient of the proof is the structure 
theorem for closed $1$-forms as it was already used in \cite{R}. 
In section 2 the connection between closed $1$-forms and
reducibility is given in two lemmas and applied to prove the theorem. 
The lemmas are proved in section 3. 

To test the quality of the estimate in the theorem we
construct examples of polynomials $f\in\Z[x,y]$ in section 4 with 
a certain reducibility behavior. Assuming the
Bouniakowsky conjecture (which will also be explained in section 
4) one gets the following result:

\begin{Proposition} Let $m,n\ge 1$ be integers. 
If the Bouniakowsy conjecture is true there are infinitely many
polynomials $f\in\Z[x,y]$ with $\deg_xf=m$ and $\deg_yf=n$ which are
absolutely irreducible over $\Q$ but reducible for a prime $p$ with 
$$p\ge H(f)^{2m}.$$
\end{Proposition}

In case $n=1$ the inequality in the theorem is $p>(2m)^m\cdot
H(f)^{2m}$. The proposition shows then that the exponent $2m$ is
best possible. 
In case $n=2$ the exponent in the theorem is $4m+1$. In
\cite{R2} it is shown that the exponent can be improved to $6$ (for
$m=2$), $6\frac{2}{3}$ (for $m=3$) and $2m$ (for $m\ge 4$). This
supports my belief that the best exponent in the theorem will
be smaller than $2mn+n-1$ if $n\ge 2$.

\section{A criterion for reducibility}

If $f(x,y)$ is a polynomial with $\deg_xf=m$ and
$\deg_yf=n$ we write $\deg f=(m,n)$. The notation $\deg f\le
(m,n)$ will mean that $\deg_x f\le m$, $\deg_yf\le n$. If it happens
that we write $\deg f\le (m,n)$ with $m<0$ or $n<0$ then $f=0$.  

The following lemmas contain our criterion for reducibility.

\begin{Lemma} Let $k$ be an arbitrary algebraically closed field and 
$f(x,y)\in k[x,y]$ a reducible polynomial with $\deg f=(m,n)$. Then
there are polynomials $r,s\in k[x,y]$ with $\deg r\le(m-1,n)$ and 
$\deg s\le(m,n-2)$ such that 
$$\frac{\partial}{\partial y}\left(\frac{r}{f}\right)=
\frac{\partial}{\partial x}\left(\frac{s}{f}\right)\quad\mbox{ and
}\quad (r,s)\ne (0,0).$$
\end{Lemma}

\begin{Lemma} Let $k$ be an arbitrary algebraically closed field 
of characteristic $0$ and 
$f(x,y)\in k[x,y]$ with $\deg f=(m,n)$ and $n\ge 1$. If 
there are polynomials $r,s\in k[x,y]$ with $\deg r\le(m-1,n)$ and 
$\deg s\le(m,n-2)$ such that 
$$\frac{\partial}{\partial y}\left(\frac{r}{f}\right)=
\frac{\partial}{\partial x}\left(\frac{s}{f}\right)\quad\mbox{ and
}\quad (r,s)\ne (0,0)$$
then $f$ is reducible. 
\end{Lemma}

The proof of the lemmas will be postponed to the next section.
We remark that the example $f=x$, $r=1$, $s=0$ shows that $n\ge 1$ 
is a necessary condition in lemma 2. 

\bigskip

We reformulate the lemmas: Let $f\in k[x,y]$ have degree $(m,n)$
and assume that $m,n\ge 1$. When do we find $r,s\in k[x,y]$ with 
$\deg r\le (m-1,n)$ and $\deg s\le (m,n-2)$ such that the equation 
\begin{equation}\frac{\partial}{\partial y}\left(\frac{r}{f}
\right)=\frac{\partial}{\partial x}\left(\frac{s}{f}\right)
\end{equation}
holds? We write 
$$f=\sum_{\substack{0\le i\le m\\0\le j\le n}}a_{ij}x^iy^j,\quad 
r=\sum_{\substack{0\le i\le m-1\\0\le j\le n}}u_{ij}x^iy^j,\quad 
s=\sum_{\substack{0\le i\le m\\0\le j\le n-2}}v_{ij}x^iy^j$$
with unknowns $u_{ij}$ $(0\le i\le m-1,0\le j\le n)$ and
$v_{ij}$ $(0\le i\le m,0\le j\le n-2)$. 
(There are $m(n+1)+(m+1)(n-1)= 
2mn+n-1$ unknowns $u_{ij}$ and $v_{ij}$ if $m,n\ge 1$.) 
Equation (1) can be written as 
$$\frac{\partial r}{\partial y}f-r\frac{\partial f}{\partial y}-
\frac{\partial s}{\partial x}f+s\frac{\partial f}{\partial x}=0.$$
We have 
$$\frac{\partial r}{\partial y}f-r\frac{\partial f}{\partial y}-
\frac{\partial s}{\partial x}f+s\frac{\partial f}{\partial x}=
\sum_{k,l}g_{kl}x^ky^l$$
with 
$$g_{kl}=\sum_{(i,j)\in A_{kl}}(-l+2j-1)a_{k-i,l-j+1}u_{ij} +
\sum_{(i,j)\in B_{kl}}(k-2i+1)a_{k-i+1,l-j}v_{ij}$$
where 
\begin{eqnarray*}
A_{kl}&=&\{(i,j):0\le k-i\le m,\quad 0\le l-j+1\le n,\quad 
0\le i\le m-1, \quad 0\le j\le n\},\\
B_{kl}&=&\{(i,j):0\le k-i+1\le m,\quad 0\le l-j\le n,\quad 
0\le i\le m, \quad 0\le j\le n-2\}.
\end{eqnarray*}
One sees that $\deg \sum g_{kl}x^ky^l\le (2m-1,2n-2)$.
Equation (1) is satisfied iff we find $u_{ij},v_{ij}\in k$ with 
$$g_{00}=\dots=g_{2m-1,2n-2}=0.$$
We can write this as a matrix equation 
$$\left(\begin{array}{c}g_{00}\\ \vdots \\ g_{2m-1,2n-2}\end{array}
\right)=
M(f)\cdot \left(\begin{array}{c}u_{00}\\ \vdots \\ u_{m-1,n} \\
v_{00} \\ \vdots \\ v_{m,n-2} \end{array}\right)=0$$
where the entries of the matrix $M(f)$ are coefficients of 
certain $g_{kl}$ with respect to $u_{ij}$ and $v_{ij}$. 

With these notations it is clear that equation (1) has a
nontrivial solution iff $M(f)$ has rank $<(2mn+n-1)$, i.e. all
$(2mn+n-1)\times (2mn+n-1)$-submatrices of $M(f)$ vanish. Now
we can reformulate the two lemmas for $f\in k[x,y]$ in terms of the
matrix $M(f)$: 
\begin{itemize}
\item If $f$ is reducible then ${\rm rank } M(f)<2mn+n-1$. 
\item If $k$ has characteristic $0$ and ${\rm rank }M(f)<2mn+n-1$ 
then $f$ is reducible. 
\end{itemize}

\bigskip

We apply this to prove the theorem: Let $f\in \Z[x,y]$ be absolutely
irreducible of degree $(m,n)$. Then the matrix $M(f)$ has rank $2mn+n-1$, i.e. there 
is a $(2mn+n-1)\times (2mn+n-1)$-submatrix $M_0$ of $M(f)$ with $\det
M_0\ne 0$. 
We will estimate $|\det M_0|$ using Hadamard's estimate for 
determinants. To do this we have to know the $L_2$-norm of the 
rows of $M_0$. A row of $M_0$ is given by the coefficients of a 
linear form $g_{kl}$ with respect to the variables $u_{ij}$ and 
$v_{ij}$. We have 
\begin{eqnarray*}
||g_{kl}||_2^2&=&\sum_{(i,j)\in A_{kl}}(-l+2j-1)^2a_{k-i,l-j+1}^2
+\sum_{(i,j)\in B_{kl}}(k-2i+1)^2a_{k-i+1,l-j}^2\le \\
&\le&(\sum_{(i,j)\in A_{kl}}(-l+2j-1)^2 +\sum_{(i,j)\in B_{kl}}
(k-2i+1)^2)\cdot H(f)^2.
\end{eqnarray*}
If $(i,j)\in A_{kl}$ then $0\le l-j+1\le n$ and $0\le j\le n$ so that 
$-n\le -(l-j+1)+j\le n$ and $(-l+2j-1)^2\le n^2$. Furthermore $\#
A_{kl}\le m(n+1)$.

If $(i,j)\in B_{kl}$ then $0\le k-i+1\le m$ and $0\le i\le m$ so that 
$-m\le (k-i+1)-i\le m$ and $(k-2i+1)^2\le m^2$. Furthermore $\#
B_{kl}\le (m+1)(n-1)$. 

This implies 
$$||g_{kl}||_2^2 \le (n^2\cdot \#A_{kl}+m^2\cdot \#B_{kl})\cdot H(f)^2\le [m(n+1)n^2+(m+1)(n-1)m^2]\cdot H(f)^2$$
so that the $L_2$-norm of a row of $M_0$ is $\le
\sqrt{[m(n+1)n^2+(m+1)(n-1)m^2]\cdot H(f)^2}$ and therefore
using Hadamard 
\begin{eqnarray*}
|\det M_0|&\le&\sqrt{[m(n+1)n^2+(m+1)(n-1)m^2]\cdot H(f)^2}^{2mn+n-1}
\\ &=&[m(n+1)n^2+(m+1)(n-1)m^2]^{mn+\frac{n-1}{2}}\cdot H(f)^{2mn+n-1}.
\end{eqnarray*}
Now if $p$ is any prime with 
\begin{equation*}
p>[m(n+1)n^2+(m+1)(n-1)m^2]^{mn+\frac{n-1}{2}}\cdot H(f)^{2mn+n-1}
\end{equation*}
then $0<|\det M_0|<p$ which implies that $\det M_0\not\equiv 0\bmod p$
so that $M(f)$ considered as a matrix over $\F_p$ has rank $2mn+n-1$ and
$f\bmod p$ is absolutely irreducible by the above criterion. 
This proves our theorem. 

\section{Proof of Lemma 1 and 2}

We start with a remark: If $k$ is an algebraically closed field and 
$g\in k[x,y]$ satisfies $\dd{g}{x}=\dd{g}{y}=0$ then $g$ is constant
in characteristic $0$ or a $p$-power in characteristic $p$. In
each case, $g$ is not irreducible. 

\bigskip

{\it Proof of Lemma 1:} Let $f\in k[x,y]$ be reducible of degree
$(m,n)$. We have to construct a nontrivial solution for the equation
$\dd{}{y}(\frac{r}{f})=\dd{}{x}(\frac{s}{f})$ with $\deg r\le (m-1,n)$
and $\deg s\le (m,n-2)$. We distinguish different cases:\\
{\sl Case I:} $f$ is squarefree. We write $f=gh$ with 
$\deg_yg=\ell$ and we can assume that $h$ is irreducible. Writing
$$g=b_0(x)+b_1(x)y+\dots+b_{\ell}(x)y^{\ell},\quad 
h=c_0(x)+c_1(x)y+\dots +c_{n-\ell}(x)y^{n-\ell}$$
gives 
\begin{eqnarray*}
\frac{\partial g}{\partial y}h&=&
b_1(x)c_0(x)+\dots+\ell b_{\ell}(x)c_{n-\ell}(x)y^{n-1},\\
g\frac{\partial h}{\partial y}&=&
b_0(x)c_1(x)+\dots+(n-\ell)b_{\ell}(x)c_{n-\ell}(x)y^{n-1}.
\end{eqnarray*}
{\sl Case I.1:} $\ell\ne 0$ in $k$. Take 
$$r=(n-\ell)\dd{g}{x}h-\ell g\dd{h}{x}\quad \mbox{ and }\quad 
s=(n-\ell)\dd{g}{y}h-\ell g\dd{h}{y}.$$
One sees at once that
$\dd{}{y}(\frac{r}{f})=\dd{}{x}(\frac{s}{f})$ holds and that by
construction $\deg r\le (m-1,n)$, $\deg s\le (m,n-2)$. 
If we had $r=s=0$ then $h$ would divide $\dd{h}{x}$ and 
$\dd{h}{y}$ which would imply $\dd{h}{x}=\dd{h}{y}=0$, 
contradicting the irreducibility of $h$. Therefore $(r,s)\ne (0,0)$ 
and we are done. \\
{\sl Case I.2:} $\ell=0$ in $k$. Then $\deg_y \dd{g}{y}h\le n-2$. 
Take $$r=\dd{g}{x}h,\quad s=\dd{g}{y}h.$$
Then the equation $\dd{}{y}(\frac{r}{f})=\dd{}{x}(\frac{s}{f})$
is satisfied with $\deg r\le (m-1,n)$ and $\deg s\le (m,n-2)$. Also
$(r,s)\ne (0,0)$ else $g$ would be a $p$-power contradicting the
fact that $f$ is supposed to be squarefree. \\
{\sl Case II:} $f$ is not squarefree. We write $f=g^2h$ and we
can assume that $g$ is irreducible. Take 
$$r=h\dd{g}{x}\quad\mbox{ and }\quad s=h\dd{g}{y}.$$
Then $(r,s)\ne (0,0)$ because 
$g$ is irreducible and 
$$\frac{r}{f}=\frac{1}{g^2}\dd{g}{x}=\dd{}{x}(-\frac{1}{g}),\quad 
\frac{s}{f}=\frac{1}{g^2}\dd{g}{y}=\dd{}{y}(-\frac{1}{g})$$
shows that $\dd{}{y}(\frac{r}{f})=\dd{}{x}(\frac{s}{f})$
holds. It is clear that $\deg r\le (m-1,n)$ and $\deg s\le
(m,n-2)$. 
\qed

{\it Proof of Lemma 2:} Suppose that $k$ is algebraically closed of
characteristic $0$, $f\in k[x,y]$ is irreducible with $\deg
f=(m,n)$ and 
$$\frac{\partial }{\partial y}\left(\frac{r}{f}\right)=
\frac{\partial }{\partial x}\left(\frac{s}{f}\right)$$
with $\deg r\le(m-1,n)$, $\deg s\le(m,n-2)$ and $(r,s)\ne (0,0)$. 
The equation implies that 
$$\omega=\frac{r}{f}dx+\frac{s}{f}dy$$
is a nontrivial closed differential form. Now the structure theorem
for closed $1$-forms (cf. \cite[Satz 2, p.172]{R}) says that 
$\omega$ has the form 
$$\omega=\sum_{i=1}^u\lambda_i\frac{dp_i}{p_i}+d(\frac{g}{q_1^{e_1}\dots
q_v^{e_v}})$$
where $p_i,q_j\in k[x,y]$ are irreducible, $g\in k[x,y]$,
$\lambda_i\in k$, $e_j\ge 0$, $p_1,\dots,p_u$ are pairwise prime,
$q_1,\dots,q_v,g$ are pairwise prime. Comparing the coefficients of
$dx$ and $dy$ gives 
\begin{eqnarray*}
\frac{r}{f}&=&\frac{\lambda_1\dd{p_1}{x}}{p_1}+\dots+\frac{\lambda_r
\dd{p_u}{x}}{p_u}+\frac{\dd{g}{x}}{q_1^{e_1}\dots q_v^{e_v}}
-\frac{e_1g\dd{q_1}{x}}{q_1^{e_1+1}q_2^{e_2}\dots q_v^{e_v}}-\dots
-\frac{e_vg\dd{q_v}{x}}{q_1^{e_1}\dots q_{s-1}^{e_{s-1}}q_v^{e_v+1}}\\
\frac{s}{f}&=&\frac{\lambda_1\dd{p_1}{y}}{p_1}+\dots+\frac{\lambda_u
\dd{p_u}{y}}{p_u}+\frac{\dd{g}{y}}{q_1^{e_1}\dots q_v^{e_v}}
-\frac{e_1g\dd{q_1}{y}}{q_1^{e_1+1}q_2^{e_2}\dots q_v^{e_v}}-\dots
-\frac{e_vg\dd{q_v}{y}}{q_1^{e_1}\dots q_{s-1}^{e_{s-1}}q_v^{e_v+1}}
\end{eqnarray*}
$k[x,y]$ is factorial and therefore we have for each $p_i$ and $q_j$ a
valuation $v_{p_i}$ and $v_{q_j}$. \\
If $g\ne 0$ and $e_j\ge 1$ for some $j$ we would get
$v_{q_j}(\frac{r}{f})=-e_j-1\le -2$ or $v_{q_j}(\frac{s}{f})=-e_j 
-1\le -2$ as $(\dd{q_j}{x},\dd{q_j}{y})\ne (0,0)$, a contradiction
to the irreducibility of $f$. Therefore we can assume
$e_1=\dots=e_v=0$. \\
If $\lambda_i\ne 0$ and $p_i$ is prime to $f$ then
$(\dd{p_i}{x},\dd{p_i}{y})\ne (0,0)$ would imply
$v_{p_i}(\frac{r}{f})=-1$ or $v_{p_i}(\frac{s}{f})=-1$, a
contradiction. We can write now 
$$\omega=\lambda\frac{df}{f}+dg$$
with $\lambda\in k$ which gives 
$$r=\lambda\dd{f}{x}+f\dd{g}{x}\quad\mbox{ and }\quad 
s=\lambda\dd{f}{y}+f\dd{g}{y}.$$
If $\dd{g}{x}\ne 0$ then $r$ would have degree $\ge m$ in $x$, a
contradiction, if $\dd{g}{y}\ne 0$ then $s$ would have degree $\ge n$
in $y$, a contradiction. Therefore we get 
$$r=\lambda\dd{f}{x}\quad\mbox{ and }\quad s=\lambda\dd{f}{y}$$
with $\lambda\ne 0$. As $n\ge 1$ we can write 
$f=a_0(x)+\dots+a_n(x)y^n$ with $a_n(x)\ne 0$ and get 
$\dd{f}{y}=a_1(x)+\dots+na_n(x)y^{n-1}$
which shows that $s$ has degree $n-1$ in $y$, a
contradiction. Therefore $f$ can not be irreducible. 
This proves the lemma. \qed

\section{Examples} 

In the following lemma families of polynomials are constructed
with an explicit reducibility condition. 

\begin{Lemma} 
\begin{enumerate}
\item Let $k$ be an algebraically closed field of characteristic 
$\ne 2$, $m,n\ge 1$ integers and $t\in k$. The polynomial 
$f_t(x,y)=(tx^m-2x+2)+(x^m-t)y^n\in k[x,y]$
is reducible if and only if $(t^2+2)^m-2^mt=0$.
In this case the factor $x-\frac{t^2+2}{2}$ splits off. 
\item The polynomial $g_m(t)=(t^2+2)^m-2^mt\in \Z[t]$
is irreducible over $\Q$ and $\gcd\{g_m(\ell):\ell\in\N\}=1$.
\end{enumerate}
\end{Lemma}

{\it Proof:} 
\begin{enumerate}
\item Suppose first that $tx^m-2x+2$ and $x^m-t$ are relatively prime
and $f_t$ is reducible. Then $f_t$ is reducible as a
polynomial in $y$ with coefficients in $k(x)$ and therefore 
$\frac{-tx^m+2x-2}{x^m-t}$ is a nontrivial power in $k(x)$. 
Then $-tx^m+2x-2$ and $x^m-t$ have to be nontrivial powers in $k[x]$ 
and therefore inseparable. But $x^m-t$ is inseparable only if
$m=0$ or $t=0$ in $k$ and for both cases $-tx^m+2x-2$ is
separable. So this case can not happen. \\
If $tx^m-2x+2$ and $x^m-t$ have a common factor
$x-u$ for some $u\in k$ then $f_t$ is clearly reducible. 
This happens iff
$tu^m-2u+2=u^m-t=0$ which is equivalent to 
$u=\frac{t^2+2}{2}$ and $(t^2+2)^m-2^mt=0$ which proves part 1
of the lemma. 
\item Let $\alpha\in\overline{\Q}$ be any root of $g_m$ over $\Q$, 
i.e.  $\alpha=(\frac{\alpha^2+2}{2})^m$. Define
$\beta=\frac{\alpha^2+2}{2}\in\Q(\alpha)$. Then 
$\alpha=\beta^m\in\Q(\beta)$ and therefore $\Q(\alpha)=\Q(\beta)$.
Finally 
$0=\alpha^2+2-2\beta=\beta^{2m}-2\beta+2$
shows that $\beta$ is a root of the irreducible Eisenstein polynomial
$t^{2m}-2t+2$, which implies that $\Q(\alpha)=\Q(\beta)$ has degree
$2m$ over $\Q$. Therefore $g_m=(t^2+2)^m-2^mt$ is irreducible over 
$\Q$. From $g_m(0)=2^m$ and $g_m(1)\equiv 1\bmod 2$ one sees
that $\gcd\{g_m(\ell):\ell\in\N\}=1$. \qed
\end{enumerate}

To construct infinitely many examples with the right reduction
behavior we use the very plausible Bouniakowsky conjecture which was 
generalized by Schinzel as hypothesis H (cf. \cite{B},\cite{S}): 

\begin{Conjecture}[Bouniakowsky] If $g(t)\in\Z[t]$ is
irreducible and $N=\gcd\{g(\ell):\ell\in\N\}$ then there are
infinitely many $\ell\in\N$ such that $\frac{1}{N}|g(\ell)|$
is a prime. 
\end{Conjecture}

Now we prove our proposition of section 1. We use the
notations and results of the previous lemma. 
Let $m,n\ge 1$ be integers and take 
$$f_{\ell}(x,y)=(\ell x^m-2x+2)+(x^m-\ell)y^n\in \Z[x,y]$$
with $\ell\in\Z,\ell\ge 2$. Then $H(f_{\ell})=\ell$.
As $g_m(\ell)\ne 0$ in $\Q$ the polynomial $f_{\ell}$ is absolutely
irreducible over $\Q$. If $p_{\ell}=g_m(\ell)$ is a prime,
then $g_m(\ell)\equiv 0\bmod p_{\ell}$ and 
$f_{\ell}\bmod p_{\ell}$ is reducible and 
$$p_{\ell}=g_m(\ell)\ge \ell^{2m}=H(f_{\ell})^{2m}.$$
Now the Bouniakowsky conjecture says that there are infinitely
many $\ell$ such that $g_m(\ell)$ is prime. This proves the
proposition.


\begin{thebibliography}{99}
\bibitem[B]{B} V. Bouniakowsky, Nouveaux th\'eor\`emes
relatifs \`a la distinction des nombres premiers et \`a la
d\'ecomposition des entiers en facteurs, M\'em. Acad. Sc. St.
P\'etersbourg (6), Sc. Math. Phys. {\bf 6} (1857), 305--329.
\bibitem[R1]{R} W. Ruppert, Reduzibilit\"at ebener Kurven,
Journal f\"ur die reine und angewandte Mathematik {\bf
369} (1986), 167--191.
\bibitem[R2]{R2} W. Ruppert, Reducibility of polynomials
$a_0(x)+a_1(x)y+a_2(x)y^2$ modulo $p$, to appear in Archiv der
Mathematik. 
\bibitem[S]{S} A. Schinzel, W. Sierp\'inski, Sur certaines
hypoth\`eses concernant les nombres premiers, Acta Arith. {\bf
4} (1958), 185--208; corrig. Acta Arith. {\bf 5} (1959), 259.
\bibitem[Z]{Z} U. Zannier, On the reduction modulo
$p$ of an absolutely irreducible polynomial $f(x,y)$, Arch.
Math. {\bf 68} (1997), 129--138.
\end{thebibliography}
\end{document}